\documentclass[11pt]{article}
\usepackage{amssymb,amsmath}
\usepackage{graphicx,units,mathrsfs,bm}
\usepackage{subfigure}
\usepackage[sort&compress,numbers]{natbib} 
\usepackage{hyperref,xcolor}
\usepackage[paperwidth=8.5in, paperheight=11in, portrait, top=1in, bottom=1.in, left=1.15in, right=1.15in]{geometry}
\usepackage{tikz}
\usetikzlibrary{calc,decorations.markings,positioning}
\allowdisplaybreaks
\newcommand\bt{\raise 2pt \hbox{$\bigtriangledown$}\hskip 1.5pt}

\newmuskip\pFqskip
\mathchardef\pFcomma=\mathcode`, 

\numberwithin{equation}{section}

\newtheorem{thm}{Theorem}[section]

\newtheorem{lemma}{Lemma}[section]
\newtheorem{coro}{Corollary}[section]

\begin{document}

\title{\bf 

Free boson realization of the Dunkl intertwining operator in one dimension
}
\author{

Luc Vinet\textsuperscript{$1,2$}\footnote{
E-mail: vinet@CRM.UMontreal.CA}~, 
Alexei Zhedanov\textsuperscript{$3$}\footnote{
E-mail: zhedanov@yahoo.com} \\[.5em]
\textsuperscript{$1$}\small~Centre de Recherches Math\'ematiques, 
Universit\'e de Montr\'eal, \\
\small~P.O. Box 6128, Centre-ville Station, Montr\'eal (Qu\'ebec), \\
\textsuperscript{$2$}\small~
Institut de valorisation des données (IVADO),\\
\small~Montréal (Québec), H2S 3H1, Canada
H3C 3J7, Canada.\\[.9em]

\textsuperscript{$3$}\small~School of Mathematics, 
Renmin University of China, Beijing, 100872, China
}
\date{\today}

\maketitle
\begin{center}
  \textit{Dedicated to the memory of Ji\v{r}i Patera}\\[.9em]  
\end{center}

\hrule
\begin{abstract}\noindent 
The operator that intertwines between the $\mathbb{Z}_2$ - Dunkl operator and the derivative is shown to have a realization in terms of the oscillator operators in one dimension. This observation rests on the fact that the Dunkl intertwining operator maps the Hermite polynomials on the generalized Hermite polynomials.
\end{abstract}

\bigskip

\hrule

\section{Introduction}
Dunkl operators \cite{dunkl1989differential} are differential-difference operators associated to reflection groups. One of their key features is that they form a commutative set. These operators are proving of high importance in many areas and in particular in the study of orthogonal polynomials \cite{dunkl2014orthogonal} and of quantum many-body systems of the Calogero-Sutherland type (see for example \cite{lapointe1996exact}). Reviews of Dunkl theory can be found in \cite{anker2017introduction}, \cite{dunkl2008reflection}, \cite{heckman1997dunkl}, \cite{rosler2003dunkl}. In its simplest expression in one dimension, where the reflection group is $\mathbb{Z}_2$, the Dunkl operator $D_{\mu}$ provides a one-parameter deformation of the ordinary derivative $\frac{d}{dx}$:
\begin{equation}
    D_{\mu} = \frac{d}{dx} + \frac{\mu}{x} (1 - R), \quad \mu \in \mathbb{R},
\end{equation}
with $R$ the reflection operator: $Rf(x)=f(-x)$. 

A key element in this framework is the Dunkl intertwining operator $V_{\mu}$ which as it name indicates intertwines between the Dunkl derivative and the ordinary one:
\begin{equation}
    D_{\mu} V_{\mu} = V_{\mu} \frac{d}{dx}.
\end{equation}
This operator $V_{\mu}$ has been introduced in \cite{dunkl1991integral}. In the $\mathbb{Z}_2$ case that we shall be considering in this note, it is readily checked that  $V_{\mu}$ can be defined by the following action in the monomial basis:
\begin{equation}
    V_{\mu} \big(x^{2n + \epsilon}\big) = \frac{\big(\frac{1}{2}\big)_{n+\epsilon}}{(\mu +\frac{1}{2})_{n+\epsilon}}\; x^{2n + \epsilon}, \qquad n=0, \dots, \quad \epsilon = 0, 1, \label{vmumon}
\end{equation}
with $(a)_n= \frac{\Gamma(a+n)}{\Gamma[n]}$ the standard Pochhammer symbol. With the help of beta integrals, it is shown \cite{dunkl1991integral} (see also the reviews \cite{dunkl2008reflection}, \cite{rosler2003dunkl}) that this amounts to the following representation:
\begin{equation}
     V_{\mu} f(x) = 2^{-2\mu}\frac{\Gamma(2\mu +1)}{\Gamma(\mu)\Gamma(\mu +1)} \int_{-1}^1f(xt)(1-t)^{\mu -1}(1+t)^{\mu} dt. \label{A_1}
\end{equation}
In the general theory, the intertwining operators are of significant interest in particular because they provide through the Dunkl kernel the simultaneous eigenfunctions of the commuting Dunkl operators. This motivates the challenging search for explicit expressions that extend \eqref{A_1} to other reflection groups. Dunkl himself has obtained formulas for the reflection groups $A_2$ \cite{dunkl1995intertwining} and $B_2$ \cite{dunkl2007intertwining} and this remains a topic of intensive investigations; see for examples the recent reports \cite{rosler2019positive}, \cite{xu2020intertwining}, \cite{de2021dunkl} and the background references therein for an overview of the state of the art.

The free field formalism makes use of bosonic or fermionic operators that obey Wick's rule to construct representations for algebras arising in various models. This typically proves fundamentally instructive and generally quite practical from the computational viewpoint. The literature on this is vast, to pick a reference we might cite a paper \cite{morozov1993q} co-authored by one of us that provides a models of the q-hypergeometric functions in this picture. 

The purpose of the present note is to offer a free field realization of the Dunkl intertwiner operator $V_{\mu}$ in one dimension using observations made by Odake in \cite{odake2020free} on the connection between the wave functions of the harmonic and singular oscillators. The analysis is rooted in the fact that the Hermite polynomials are mapped \cite{dunkl1991integral}, \cite{cheikh2007characterization}, \cite{vinet2011missing} onto the generalized Hermite polynomials \cite{szeg1939orthogonal}, \cite{chihara2011introduction}, \cite{rosenblum1994generalized} by the intertwining operator $V_{\mu}$.

The paper is organized as follows. The generalized Hermite polynomials are introduced next through the eigenfunctions of the Dunkl oscillator. The key finding is provided in Section \ref{3} where the action of the intertwiner $V_{\mu}$ in the standard Hermite polynomial basis is explicitly obtained. This is then translated in terms of the oscillator annihilation operator in Section \ref{4} to obtain the desired free boson realization of $V_{\mu}$. Remarks on possible generalizations are offered in guise of conclusion.

\section{The Dunkl oscillator and the generalized Hermite polynomials}
The one-dimensional Dunkl oscillator is governed by the Hamiltonian \cite{genest2013dunkl}
\begin{equation}
    H_{\mu} = \frac{1}{2}\big( -D_{\mu}^2 + x^2\big) = \frac{1}{2} \Big[-\frac{d^2}{dx^2} - \frac{2\mu}{x}\; \frac{d}{dx} + \frac{\mu}{x^2}(1 - R) + x^2 \Big].
\end{equation}
Its eigenfunctions $\psi^{(\mu)}_{2n+\epsilon}(x)$ with eigenvalues $E= 2n + \epsilon + \mu +\frac{1}{2}$ are given by
\begin{equation}
    \psi^{(\mu)}_{2n+\epsilon}(x)= \kappa_n e^{-\frac{1}{2}x^2} H_{2n+\epsilon}^{\mu}(x), \quad n=0, 1, \dots
\end{equation}
with $\kappa_n$ a normalization constant. The functions $H_{2n+\epsilon}^{\mu}(x)$ are the generalized Hermite polynomials which can be defined in terms of the Laguerre polynomials \cite{koekoek2010hypergeometric}:
\begin{equation}
   L_n^{\alpha} (x) = \frac{(\alpha + 1)_n}{n!} {_1}F_1 \left( {-n \atop \alpha + 1} ;x\right) = \frac{(\alpha + 1)_n}{n!}\:\sum_{k=0}^n \frac{(-n)_k}{(\alpha + 1)_k} \frac{x^k}{k!}; \label{lagexp}
\end{equation}
they read \cite{cheikh2007characterization}, \cite{rosenblum1994generalized} 
\begin{equation}
    H_{2n+\epsilon}^{\mu}(x) = \frac{(-1)^n (2n + \epsilon)!}{(\mu + \frac{1}{2})_{n+\epsilon}} x^{\epsilon} L_m^{\mu - \frac{1}{2} + \epsilon}(x^2). \label{genHer}
\end{equation}
When $\mu = 0$, the Dunkl oscillator manifestly reduces to the ordinary harmonic oscillator whose Hamiltonian as is well known can be written in the form $H_{\mu=0} \equiv N +\frac{1}{2}$ = with $N$ and the annihilation and creation operators given by
\begin{equation}
   N=a^{\dagger}a,\qquad a=\frac{1}{\sqrt{2}} (x + \frac{d}{dx}), \qquad a^\dagger = \frac{1}{\sqrt{2}} (x - \frac{d}{dx}).
\end{equation}
Moreover the generalized Hermite polynomials reduce to the ordinary Hermite polynomials $H_n(x)$. From \eqref{genHer} and the explicit expression of the Laguerre polynomials, one readily obtains the following formula that will prove useful for the Hermite polynomials:
\begin{equation}
    H_{2n+\epsilon}(x) = \sum_{k=0}^n C_{nk}^{(\epsilon)} x^{2k+\epsilon}, \qquad n=0, 1, \dots \label{Herexpl}
\end{equation}
with the coefficients $C_{nk}^{(\epsilon)}$ given by
\begin{equation}
    C_{nk}^{(\epsilon)}=\frac{(-1)^{n-k}\; 2^{2k+\epsilon}\;(2n+\epsilon)!}{(2k+\epsilon)!\;(n-k)!}.\label{Hercoef}
\end{equation}
The Hermite polynomials are well known to possess the Appell property
\begin{equation}
    \frac{d}{dx} H_n(x) = 2n H_{n-1}(x) \label{app}
\end{equation}
which is readily derived from the explicit expression given in \eqref{Herexpl}, \eqref{Hercoef}. They also satisfy \cite{koekoek2010hypergeometric} the eigenvalue equation
\begin{equation}
    \hat{N}H_n(x)=\frac{1}{2}\big{(}-\frac{d^2}{dx^2}+2x\big{)}H_n(x)=n H_n(x).\label{Ham}
\end{equation}
Note that $\frac{d}{dx}$ and $\hat{N}$ are related to the standard bosonic operators by a simple conjugation:
\begin{equation}
    \hat{a}=\frac{1}{\sqrt{2}}\dfrac{d}{dx}=e^{\frac{1}{2}x^2}ae^{-\frac{1}{2}x^2} \qquad 
    \hat{N}=e^{\frac{1}{2}x^2}a^{\dagger}ae^{-\frac{1}{2}x^2}. \label{bosop}
\end{equation}
\section{The Dunkl intertwiner and the Hermite polynomials} \label{3}

As an essential step towards bosonizing the Dunkl intertwiner, we shall obtain its action in the basis of the Hermite polynomials. This will in fact provide an expansion of the generalized Hermite polynomials $H^{\mu}_n(x)$ in terms of the ordinary ones with $\mu=0$. The strategy is simple: apply the definition \eqref{vmumon} of $V_{\mu}$ to \eqref{Herexpl} and use the inversion of this last formula to revert to Hermite polynomials. To that end, let us start with the following:
\begin{lemma} \label{lem}
The relation 
\begin{equation}
    \sum^m_{n=0}\frac{(-m)_n(-n)_k}{n!}=m!\; \delta_{m,k}
\end{equation}
holds for $m$ and $k$ non-negative integers.
\end{lemma}
\textit{Proof}. First note that the sum effectively starts at $n=k$. Observe then that 
\begin{equation}
    \sum_{n=k}^m \frac{(-m)_n(-n)_k}{n!}=(-1)^k(-m)_k \sum_{k=0}^{m-k} (-1)^k {m-k \choose k}
    = (-1)^k(-m)_k (1-1)^{m-k},
\end{equation}
from where we see that
\begin{equation}
    \sum_{n=k}^m \frac{(-m)_n(-n)_k}{n!}=(-1)^m(-m)_m \;\delta_{m,k}
\end{equation}
and therefore that the lemma is true.
A first application provides the inversion of formula \eqref{Herexpl}, namely the expansion of monomials in terms of Hermite polynomials \cite{luke2014mathematical}, \cite{sanchez1998expansions}. It is given by
\begin{equation}
    x^{2k+\epsilon}=\sum_{l=0}^k D_{k\ell}^{(\epsilon)} H_{2\ell +\epsilon}, \qquad k=0, 1, \dots \label{invers}
\end{equation}
with the the coefficients $D_{k\ell}^{(\epsilon)}$ reading
\begin{equation}
   D_{k\ell}^{(\epsilon)}=\frac{(2k+\epsilon)!}{2^{2k+\epsilon}(2\ell+\epsilon)!(k-\ell)!}. \label{Dentries}
\end{equation}
This can be confirmed by showing that the matrix $D$ with entries given in \eqref{Dentries} is the inverse of the matrix $C$ made out of the coefficients provided by \eqref{Hercoef}. Indeed using
\begin{equation}
    (n-k)!=\frac{n!}{(-1)^k(-n)_k},\label{eqn0}
\end{equation}
it is readily checked that 
\begin{align}
        \sum_{k=0}^nC_{nk}D_{k\ell} &= \sum_{k=0}^n(-1)^{n-k} \frac{(2n)!}{(2\ell) !}\frac{1}{(n-k)! (k-\ell)!}\\
        &=(-1)^{n-\ell} \frac{(2n)!}{(2\ell)!} \frac{1}{n!}\sum_{k=0}^n \frac{(-n)_k (-k)_l}{k!} = \delta _{n,\ell},
\end{align}
 calling upon Lemma~\ref{lem} at the last step.
 
 In preparation for the computation of $V_{\mu}H_{2n+\epsilon}(x)$, let us record a few useful formulas.
 \begin{lemma} \label{lem2}
 The equation
 \begin{equation}
 \frac{(\frac{1}{2})_{k+\epsilon}}{(\mu +\frac{1}{2})_{k+\epsilon}} = \frac{(\frac{1}{2})_{n+\epsilon}}{(\mu +\frac{1}{2})_{n+\epsilon}} \sum_{s=0}^{n-k} \frac{1}{s!} (\mu)_s \frac{(-n+k)_s}{(-n-\epsilon + \frac{1}{2})_s} \label{id2}
 \end{equation}
is an identity for $k$ and $n$ integers such that $k\le n$.
 \end{lemma}
 \textit{Proof}. Note that
 \begin{equation}
    \sum_{s=0}^{n-k} \frac{1}{s!} (\mu)_s \frac{(-n+k)_s}{(-n-\epsilon + \frac{1}{2})_s} = {_2}F_1 \left( {-n + k, \:\mu \atop -n - \epsilon + \frac{1}{2} }\; ;1\right) = \frac{(-n-\epsilon +\frac{1}{2} - \mu)_{n-k}}{(-n-\epsilon +\frac{1}{2})_{n-k}} \label{eqn1}
 \end{equation}
 using Gauss' summation formula \cite{gasper2004basic}. Now from the following properties of the Pochhammer symbols:
 \begin{align}
     (a)_{n-k}&=\frac{(-1)^k(a)_n}{(1-n-a)_k},\\
     (-n+\frac{1}{2}-a)_n&=(-1)^n (a + \frac{1}{2})_n,\\
     \frac{(a+\epsilon)_k}{(a+\epsilon)_n}&=\frac{(a)_{k+\epsilon}}{(a)_{n+\epsilon}},
 \end{align}
we immediately have that
\begin{equation}
    \frac{(-n-\epsilon +\frac{1}{2} - \mu)_{n-k}}{(-n-\epsilon +\frac{1}{2})_{n-k}}=\frac{(\mu + \frac{1}{2})_{n+\epsilon}}{(\mu + \frac{1}{2})_{k+\epsilon}} \cdot \frac{(\frac{1}{2})_{k+\epsilon}}{(\frac{1}{2})_{n+\epsilon}},
\end{equation}
which in view of \eqref{eqn1} yields \eqref{id2}. We shall also need:
\begin{lemma} \label{lem3}
The formula
\begin{equation}
    (2n-2\ell+\epsilon)!\;(-n-\epsilon+\frac{1}{2})_{\ell}=\frac{(-1)^{\ell}\;(2n+\epsilon)!\;(n-\ell)!}{2^{2\ell}n!} \label{id3}
\end{equation}
is valid for $\ell \le n$ and $\epsilon = 0, 1$.
\end{lemma}
\textit{Proof}. Start with
\begin{equation}
    (2n-2\ell + \epsilon)! = \frac{(2n+\epsilon)!}{(-2n-\epsilon)_{\ell}}. \label{eqn3}
\end{equation}
Observe then that
\begin{equation}
    (-2n-\epsilon)_{2\ell}=2^{2\ell}\;(-n)_{\ell}\;(-n-\epsilon+\frac{1}{2})_{\ell}
\end{equation}
by verifying for instance the relation for $\epsilon=0$ and $\epsilon=1$. Combining with \eqref{eqn3} leads to \eqref{id3} bearing in mind \eqref{eqn0}.
We are now ready to present the central result.
\begin{thm} \label{Theo}
The action of the Dunkl intertwining operator $V_{\mu}$ in the basis of the Hermite polynomials $H_n(x)$ is given by
\begin{equation}
    V_{\mu} H_{2n+\epsilon}(x)= \frac{(\frac{1}{2})_{n+\epsilon}}{(\mu +\frac{1}{2})_{n+\epsilon}}\,\sum_{\ell = 0}^n (-1)^{\ell}\, 2^{2\ell}\, {n \choose \ell}\,(\mu)_{\ell}\, H_{2(n-\ell)+\epsilon}(x). \label{main}
\end{equation}
\end{thm}
\textit{Proof}. From \eqref{Herexpl}, \eqref{vmumon} and \eqref{invers} we have
\begin{equation}
    V_{\mu} H_{2n+\epsilon}(x)= \sum_{k=0}^n \frac{(\frac{1}{2})_{k+\epsilon}}{(\mu +\frac{1}{2})_{k+\epsilon}}\sum_{\ell=0}^k C_{nk}^{(\epsilon)}\,D_{k\ell}^{(\epsilon)}\, H_{2\ell +\epsilon}(x).
\end{equation}
Exchanging the sums and using the expressions for $C_{nk}^{(\epsilon)}$ and $D_{l\ell}^{(\epsilon)}$ given by \eqref{Hercoef} and \eqref{Dentries} respectively, one finds after some simplifications that
\begin{equation}
    V_{\mu} H_{2n+\epsilon}(x)=\sum_{\ell=0}^n \frac{1}{(2\ell+\epsilon)!}\,H_{2\ell + \epsilon}(x)\, \sum _{k=\ell}^n \frac{(\frac{1}{2})_{k+\epsilon}}{(\mu +\frac{1}{2})_{k+\epsilon}}\,\frac{(-1)^{n-k} (2n+\epsilon)!}{(n-k)!\,(k-\ell)!}.
\end{equation}
At this point we call upon Lemma \ref{lem2} to write
\begin{align}
    V_{\mu} H_{2n+\epsilon}(x)&=(-1)^n \frac{(\frac{1}{2})_{n+\epsilon}}{(\mu +\frac{1}{2})_{n+\epsilon}}\, (2n+\epsilon)! \sum_{\ell =0}^n \frac{1}{(2\ell+\epsilon)!}\,H_{2\ell + \epsilon}(x)\nonumber\\
    &\cdot \sum_{k=\ell}^n \frac{(-1)^k}{(n-k)!\,(k-\ell)!}\, \sum_{s=0}^{n-k}\frac{1}{s!} (\mu)_s \frac{(-n+k)_s}{(-n-\epsilon + \frac{1}{2})_s} .
\end{align}
The second sum can be made to start at $k=0$ since all the terms with $k \le \ell$ vanish. One then exchanges the sums over $k$ and $s$ to have
\begin{align}
     V_{\mu} H_{2n+\epsilon}(x)&=\frac{(-1)^n}{n!}\,\frac{(\frac{1}{2})_{n+\epsilon}}{(\mu +\frac{1}{2})_{n+\epsilon}}\, (2n+\epsilon)! \sum_{\ell =0}^n \frac{(-1)^{\ell}}{(2\ell+\epsilon)!}\,H_{2\ell + \epsilon}(x)\nonumber\\
     &\cdot \sum_{s=0}^n\frac{1}{s!}\,\frac{(\mu)_s}{(-n -\epsilon +\frac{1}{2})_s} \cdot \sum_{k=0}^{n-s}\frac{1}{k!}\,(-n)_k\,(-k)_{\ell}\,(-n+k)_s \label{intstep}
\end{align}
with the help of \eqref{eqn0}. This is where we use again Lemma \ref{lem} to find for the last sum:
\begin{equation}
    \sum_{k=0}^{n-s} \frac{1}{k!}\,(-n)_k\,(-k)_{\ell}\,(-n+k)_s = (-n)_s \sum_{k=0}^{n-s} \frac{(-k)_{\ell}\,(-n+s)_k}{k!} = (-n)_s\,(n-s)!\,\delta_{n-s,\ell}. \label{intsum}
\end{equation}
Inserting \eqref{intsum} in \eqref{intstep} yields:
\begin{align}
 V_{\mu} H_{2n+\epsilon}(x)&=\nonumber\\
 &\frac{(-1)^n}{n!}\,\frac{(\frac{1}{2})_{n+\epsilon}}{(\mu +\frac{1}{2})_{n+\epsilon}}\, (2n+\epsilon)! \sum_{\ell =0}^n\frac{(-1)^{\ell}}{(2\ell + \epsilon)!}\, \frac{(\mu)_{n-\ell}}{(n-\ell)!}\,\frac{(-n)_{n-\ell}\;\ell !}{(n-\epsilon +\frac{1}{2})_{n-\ell}} H_{2\ell +\epsilon}(x).
 \end{align}
 Now effect the change of summation index $\ell \rightarrow n-\ell$ to obtain:
 \begin{equation}
    V_{\mu} H_{2n+\epsilon}(x)= \frac{(\frac{1}{2})_{n+\epsilon}}{(\mu +\frac{1}{2})_{n+\epsilon}}\, (2n+\epsilon)! \, \sum_{\ell=0}^n \frac{(\mu)_{\ell}}{\ell !}\, \frac{1}{(2n-2\ell +\epsilon)!\,(-n-\epsilon +\frac{1}{2})_{\ell}} H_{2(n-\ell) +\epsilon}(x),
 \end{equation}
 where $(-1)^{\ell}(-n)_{\ell}(n-\ell)!$ has cancelled $n!$ as is familiar from \eqref{eqn0}.
 This is where Lemma \ref{lem3} now comes in handy and allows to rewrite the denominator under the sum and to find the result \eqref{main} of Theorem \ref{Theo} after the cancellation of the factor $(2n+\epsilon)!$.
 \begin{coro}
 Equation \eqref{main} provides the connection formula between the generalized Hermite polynomials and the standard Hermite polynomials.
 \end{coro}
 \textit{Proof}. Following Dunkl \cite{dunkl1991integral}, it is easy to check that
 \begin{equation}
     H_{2n+\epsilon}^{\mu}(x) = V_{\mu} H_{2n+\epsilon}(x),
 \end{equation}
 from where the statement is established. Indeed, from \eqref{vmumon}, \eqref{lagexp} and \eqref{genHer}, we see for instance when $\epsilon=1$ that
 \begin{align}
     V_{\mu} H_{2n+1}(x)=&(-1)^n\frac{(2n+1)!}{(\frac{1}{2})_{n+1}}\,\frac{(\frac{3}{2})_n}{n!}\,\sum_{k=0}^n\frac{(-n)_k}{(\frac{3}{2})_k}\, \frac{V_{\mu}(x^{2k+1})}{k!}\nonumber\\
     =&(-1)^n\frac{(2n+1)!}{(\frac{1}{2})}\,\frac{1}{n!}\,\sum_{k=0}^n\frac{(-n)_k}{(\frac{3}{2})_k}\, \frac{(\frac{1}{2})_{k+1}}{k!\,(\mu+\frac{1}{2})_{k+1}} \, x^{2k+1}\nonumber\\
     =&(-1)^n\frac{(2n+1)!}{(\mu+\frac{1}{2})\,n!}\,x\,\sum_{k=0}^n\frac{(-n)_k}{k!}\,\frac{1}{(\mu + \frac{3}{2})_k}\,x^{2k}\nonumber\\
     =&(-1)^n\,\frac{(2n+1)!}{(\mu +\frac{1}{2})_{n+1}} \,x\,L_n^{(\mu + \frac{1}{2})}(x^2)\,=\,  H_{2n+1}^{\mu}(x),
 \end{align}
 and the computation is even more straigtforward for $\epsilon=0$.
\section{The boson operator realization} \label{4}
Borrowing ideas presented in \cite{odake2020free}, we offer in this section a realization of the Dunkl intertwining operator $V_{\mu}$ in terms of bosonic operators. This follows from Theorem \ref{Theo} which gives the action of $V_{\mu}$ on the Hermite polynomials that provide a well known representation basis for the oscillator operators. Recall the actions \eqref{app}, \eqref{Ham} on $H_n(x)$ of the bosonic operators $\hat{a}$ and $\hat{N}$ defined in \eqref{bosop}. We shall also make use of the projector
\begin{equation}
    P=\frac{1}{2} (I-R), \qquad P^2=P
\end{equation}
with $R$ the reflection operator that appears in the Dunkl derivative. Owing to the fact that the Hermite polynomials are symmetric, we have
\begin{equation}
    P\,H_{2n+\epsilon}(x)=\epsilon \,H_{2n+\epsilon}.
\end{equation}
Consider now the operator
\begin{equation}
    \hat{b}=(\hat{N} + P  + 1)^{-1}\, \hat{a}^2
\end{equation}
which clearly will lower degrees by two. It is checked that for $\epsilon =0, 1$:
\begin{align}
  \hat{b}\, H_{2n+\epsilon}(x)&=\frac{2(2n+\epsilon)(2n+\epsilon-1)}{(2n+2 \epsilon -1)}\, H_{2n+\epsilon -2} \nonumber\\
  &=2^2 \,n\, H_{2n+ \epsilon - 2}.
\end{align}
Note that $(\hat{N} + P  + 1)^{-1}$ is defined since the operator of which it is meant to be the inverse does not have zero for eigenvalue. Iterating, we have
\begin{equation}
    \hat{b}^{\ell} \, H_{2n+\epsilon}(x) = 2^{2\ell}\, \frac{n!}{(n-\ell)!}\, H_{2(n-\ell)+\epsilon}.
\end{equation}
This allows to write eq. \eqref{main} of Theorem \ref{Theo} in the form:
\begin{equation}
    V_{\mu}\,H_{2n+\epsilon}=\frac{(\frac{1}{2})_{n+\epsilon}}{(\mu +\frac{1}{2})_{n+\epsilon}}\, \sum_{\ell = 0}^n \frac{(-1)^{\ell}}{\ell !}\, \hat{b}^{\ell}\, H_{2n+\epsilon}(x).
\end{equation}
We here observe that the sum in the formula above can be extended to infinity since
\begin{equation}
    \hat{b}^{\ell}\, H_{2n + \epsilon}(x) = 0 \qquad \text{for} \qquad \ell = n+1, n+2, \dots.
\end{equation}
We thus arrive at an expression for $V_{\mu}$ in terms of oscillator operators.
\begin{thm}
On the Hilbert space $L^2(\mathbb{R}, e^{-x^2}dx)$, the Dunkl intertwining operator $V_{\mu}$ admits the following bosonic realization
\begin{equation}
    V_{\mu}= {_1}F_{0} \left( {\mu \atop -} ;-\hat{b}\right),
\end{equation}
with ${_1}F_{0}$ referring to the usual notation of hypergeometric functions \cite{koekoek2010hypergeometric, gasper2004basic}. 
\end{thm}
\section{Conclusion}
This paper has offered observations on the Dunkl intertwining operator $V_{\mu}$ in one dimension. It has provided its action in the basis of Hermite polynomials and, as a result, identified its realization in terms of bosonic operators. Knowing that the operator $V_{\mu}$ maps the Hermite polynomials on the generalized Hermite polynomials, this therefore gave the connection formula relating the two families of polynomials.

This study has a kinship with the work performed in \cite{atakishiyeva2010lifting} aimed at constructing the higher families of the Askey scheme \cite{koekoek2010hypergeometric} from the lower ones and at adding parameters (like $\mu$ here) from the action of functions of operators identified from connection formulas.

It suggests also that more cases would deserve an analysis similar to the one performed here.
Indeed we may mention two cases of polynomial families related by the Dunkl intertwiner: the pair formed by the generalized Gegenbauer polynomials and the ordinary ones \cite{cheikh2007characterization} and as well the tandem made out of the little $-1$ Jacobi polynomials and special Jacobi polynomials \cite{vinet2011missing}. Exploring connection formulas in these instances would certainly be worthwhile.

From the standpoint of 1-D quantum mechanical systems, the Dunkl intertwiner maps the wave functions of the harmonic oscillator into those of the Dunkl oscillator. It is not difficult to see that a gauge transform of the Dunkl oscillator Hamiltonian corresponds to the singular (or radial) oscillator with a reflection dependent centrifugal term \cite{genest2013hahn}. Sticking to one parity sector, say the even one, at the price of not having a unified framework, one finds that the wave functions of the singular oscillator are obtained from those of the harmonic oscillator with an even number of excitations. This is the view taken in \cite{odake2020free}. This observation is in keeping with the fact that on the one hand, the states of the harmonic oscillator with a fixed parity support metaplectic representations of $\mathfrak{su}(1,1)$ which is then mapped onto the $\mathfrak{su}(1,1)$ irreducible representation space spanned by the states of the singular oscillator. On the other hand, the Dunkl oscillator like the harmonic oscillator exhibits $\mathfrak{osp}(1|2]$ supersymmetry as expressed by the mapping of the entire Hilbert space of both systems into one another by the Dunkl intertwiner.

Finally, looking at multivariate situations should retain attention. There have been various studies like \cite{awata1995excited} connecting integrable models such as the Calogero-Sutherland one and the Virasoro algebra (or generalizations) through bosonization. This had significant impact on the study of symmetric functions. We may then ask the question of what bearing could bosonization further have on the study of Dunkl intertwining operators in this context. We hope to have modestly instilled interest in these various questions.

\section*{Acknowledgments}
 The work of LV is supported in part by a Discovery Grant from the Natural Sciences and Engineering Research Council (NSERC) of Canada. AZ who is funded by the National Foundation of China (Grant No.11771015) gratefully acknowledges the hospitality of the CRM over an extended period and the award of a Simons CRM professorship.
\bibliographystyle{unsrt} 
\bibliography{ref_ff.bib}

\begin{thebibliography}{10}

\bibitem{dunkl1989differential}
C.~F. Dunkl.
\newblock {Differential-difference operators associated to reflection groups}.
\newblock {\em Transactions of the American Mathematical Society},
  311(1):167--183, 1989.

\bibitem{dunkl2014orthogonal}
C.~F. Dunkl and Y.~Xu.
\newblock {\em {Orthogonal polynomials of several variables}}.
\newblock Number 155. Cambridge University Press, 2014.

\bibitem{lapointe1996exact}
L.~Lapointe and L.~Vinet.
\newblock {Exact operator solution of the Calogero-Sutherland model}.
\newblock {\em Communications in mathematical physics}, 178(2):425--452, 1996.

\bibitem{anker2017introduction}
J.-Ph. Anker.
\newblock {An introduction to Dunkl theory and its analytic aspects}.
\newblock In {\em Analytic, Algebraic and Geometric Aspects of Differential
  Equations}, pages 3--58. Springer, 2017.

\bibitem{dunkl2008reflection}
C.~F. Dunkl.
\newblock {Reflection groups in analysis and applications}.
\newblock {\em Japan. J. Math}, 3:215--246, 2008.

\bibitem{heckman1997dunkl}
G.J. Heckman.
\newblock Dunkl operators.
\newblock {\em Ast\'erisque}, pages 223--223, 1997.

\bibitem{rosler2003dunkl}
M.~R{\"o}sler.
\newblock {Dunkl operators: theory and applications}.
\newblock In {\em Orthogonal polynomials and special functions}, pages 93--135.
  Springer, 2003.

\bibitem{dunkl1991integral}
C.~F. Dunkl.
\newblock {Integral kernels with reflection group invariance}.
\newblock {\em Canadian Journal of Mathematics}, 43(6):1213--1227, 1991.

\bibitem{dunkl1995intertwining}
C.~F. Dunkl.
\newblock {Intertwining operators associated to the group $S_3$}.
\newblock {\em Transactions of the American Mathematical Society},
  347(9):3347--3374, 1995.

\bibitem{dunkl2007intertwining}
C.~F. Dunkl.
\newblock {An intertwining operator for the group $B_2$}.
\newblock {\em Glasgow Mathematical Journal}, 49(2):291--319, 2007.

\bibitem{rosler2019positive}
M.~R{\"o}sler and M.~Voit.
\newblock {Positive intertwiners for Bessel functions of type B}.
\newblock {\em arXiv:1912.12711}, 2019.

\bibitem{xu2020intertwining}
Y.~Xu.
\newblock Intertwining operator associated to symmetric groups and summability
  on the unit sphere.
\newblock {\em arXiv:2004.08727}, 2020.

\bibitem{de2021dunkl}
H.~De~Bie and P.~Lian.
\newblock {The Dunkl kernel and intertwining operator for dihedral groups}.
\newblock {\em Journal of Functional Analysis}, 280(7):108932, 2021.

\bibitem{morozov1993q}
A.~Morozov and L.~Vinet.
\newblock {q-Hypergeometric functions in the formalism of free fields}.
\newblock {\em Modern Physics Letters A}, 8(30):2891--2902, 1993.

\bibitem{odake2020free}
S.~Odake.
\newblock {Free Oscillator Realization of the Laguerre Polynomial}.
\newblock {\em arXiv:2008.10756}, 2020.

\bibitem{cheikh2007characterization}
Y.~Ben Cheikh and M.~Gaied.
\newblock {Characterization of the Dunkl-classical symmetric orthogonal
  polynomials}.
\newblock {\em Applied mathematics and computation}, 187(1):105--114, 2007.

\bibitem{vinet2011missing}
L.~Vinet and A.~Zhedanov.
\newblock {A ‘missing’ family of classical orthogonal polynomials}.
\newblock {\em Journal of Physics A: Mathematical and Theoretical},
  44(8):085201, 2011.

\bibitem{szeg1939orthogonal}
G.~Szeg\H{o}.
\newblock {\em {Orthogonal polynomials}}, volume~23.
\newblock American Mathematical Soc., 1939.

\bibitem{chihara2011introduction}
T.~S. Chihara.
\newblock {\em An introduction to orthogonal polynomials}.
\newblock Courier Corporation, 2011.

\bibitem{rosenblum1994generalized}
M.~Rosenblum.
\newblock {Generalized Hermite polynomials and the Bose-like oscillator
  calculus}.
\newblock In {\em Non selfadjoint operators and related topics}, pages
  369--396. Springer, 1994.

\bibitem{genest2013dunkl}
V.~X. Genest, M.~E.~H. Ismail, L.~Vinet, and A.~Zhedanov.
\newblock {The Dunkl oscillator in the plane: I. Superintegrability, separated
  wavefunctions and overlap coefficients}.
\newblock {\em Journal of Physics A: Mathematical and Theoretical},
  46(14):145201, 2013.

\bibitem{koekoek2010hypergeometric}
R.~Koekoek, P.~A. Lesky, and R.~F. Swarttouw.
\newblock {\em Hypergeometric orthogonal polynomials and their q-analogues}.
\newblock Springer Science \& Business Media, 2010.

\bibitem{luke2014mathematical}
Y.~L. Luke.
\newblock {\em Mathematical functions and their approximations}.
\newblock Academic Press, 2014.

\bibitem{sanchez1998expansions}
J.~S{\'a}nchez-Ruiz and J.~S. Dehesa.
\newblock Expansions in series of orthogonal hypergeometric polynomials.
\newblock {\em Journal of computational and applied mathematics},
  89(1):155--170, 1998.

\bibitem{gasper2004basic}
G.~Gasper and M.~Rahman.
\newblock {\em Basic hypergeometric series}, volume~96.
\newblock Cambridge university press, 2004.

\bibitem{atakishiyeva2010lifting}
M.~Atakishiyeva and N.~Atakishiyev.
\newblock {On lifting q-difference operators in the Askey scheme of basic
  hypergeometric polynomials}.
\newblock {\em Journal of Physics A: Mathematical and Theoretical},
  43(14):145201, 2010.

\bibitem{genest2013hahn}
V.~X. Genest, J.-M. Lemay, L.~Vinet, and A.~Zhedanov.
\newblock {The Hahn superalgebra and supersymmetric Dunkl oscillator models}.
\newblock {\em Journal of Physics A: Mathematical and Theoretical},
  46(50):505204, 2013.

\bibitem{awata1995excited}
H.~Awata, Y.~Matsuo, S.~Odake, and J.~Shiraishi.
\newblock {Excited states of the Calogero-Sutherland model and singular vectors
  of the $W_N$ algebra}.
\newblock {\em Nuclear Physics B}, 449(1-2):347--374, 1995.

\end{thebibliography}

\end{document}